# Matemáticas, espacios públicos e integración vecinal. El caso de Cuernavaca (México)

Igor Barahona[a*], Lucía López de Medrano[b], Barbara Martínez Moreno[c], Beatríz Limón Gutiérrez[d]

[a] Instituto de Matemáticas, Unidad Cuernavaca - Cátedras CONACyT
[b] Instituto de Matemáticas, Unidad Cuernavaca, Universidad Nacional Autónoma de México.
[c] Secretaría de Cultura del Estado de Morelos.
[d] Universidad Autónoma del Estado de Morelos.

Resumen

Investigamos el impacto de las matemáticas en el mejoramiento de la integración vecinal y la percepción de seguridad. Tomamos como caso de estudio la plaza principal colonia Chamilpa en Cuernavaca, México. La ciudad se caracteriza por tener una precaria infraestructura de recreación y ocio. La colonia se encuentra entre las de mayor marginación de la ciudad. Los datos se recogieron a través de una encuesta aplicada a los asistentes del festival ARTEMAT. Los resultados proporcionan evidencia empírica para soportar que la realización de actividades matemáticas en espacios públicos, aumenta la percepción de seguridad y mejora la cohesión social.

Palabras clave:

- *Arte*
- *Pensamiento matemático*
- *Espacios públicos*
- *Cohesión social*
- *Creatividad*

ABTRACT

We investigate the impact of mathematics on improving neighbourhood integration and perception of security. The main square of Chamilpa colony in Cuernavaca, México is take as study case. This city is featured by its precarious recreational and leisurial infrastructure. Chamilpa is among the top of social outcast levels in the city. Data was collected through a questionnaire applied among attendees of the ARTEMAT festival. Results provide empirical evidence for supporting that performing mathematics activities on public spaces, increases the perception of security and improve the social cohesion.

KEY WORDS:

- *Art*
- *Mathematical thinking*
- *Public spaces*
- *Social cohesion*
- *Creativity*

* Autor de correspondencia. Instituto de Matemáticas, Sede Cuernavaca, Universidad Nacional Autónoma de México, Av. Universidad s/n. Col. Lomas Chamilpa Código Postal 62210, Cuernavaca, Morelos. México. Correo electrónico: igor@im.unam.mx

# 1. INTRODUCCIÓN

En su tesis, Vanessa Reyes Calderón describe un acontecimiento que ilustra el nivel de degradación en el que se encuentran los espacios públicos de las principales áreas urbanas de México. –"¡*Ahí mire! ¡Junto a la barda esta, ahí es donde venden*!" – Comenta un vecino de la colonia Barrios Sierra en la delegación Magdalena Contreras de la Ciudad de México. Lo anterior durante entrevista realizada para recopilar datos relacionados con aspectos de seguridad y mantenimiento de la Unidad deportiva "Las Torres". La venta de drogas es un hecho al que se han acostumbrado los habitantes de esta colonia. El testimonio de uno de los entrevistados por Reyes (2014) proporciona elementos adicionales de análisis y reflexión. –"*Uno mejor ni se mete, hacemos como si no viéramos nada, porque si llamamos a la patrulla, hasta ellos mismos nos acusan con los delincuentes*".

De acuerdo con Ramírez (2015), la subordinación de lo público a lo privado, así como el predomino de la propiedad privada sobre el interés común, ha cambiado radicalmente la forma en que la sociedad concibe los espacios públicos. Los modelos económicos preponderantes, basados en el libre mercado, y los cuales privilegian la propiedad privada, cambian radicalmente el sentido de los espacios públicos. Bajo este paradigma es común observar mayor preocupación por mantener limpia la cochera de casa, que el jardín público que se encuentra en la colonia. Las manifestaciones de la subordinación de lo público a lo privado son variadas: espacios públicos inseguros, sucios y descuidados, sólo por mencionar algunas. Aunado a lo anterior, el mantenimiento y mejora de los espacios públicos han estado ausentes de las políticas públicas urbanas en México (Ramírez, 2014). Durante las últimas tres décadas, las políticas de desarrollo urbano han concedido a los espacios públicos como elementos de equipamiento (en fraccionamientos y unidades habitacionales), áreas verdes, de vía pública y residuales. Autores como Ramírez (2015), Ramírez (2014) y Bernard y Rowles (2013), afirman que uno de los principales problemas con los espacios públicos en México es que, desde su concepción y diseño, se encuentran disociados de la cohesión social e integración vecinal.

Recientemente, grupos de vecinos, líderes sociales y organizaciones civiles han reclamado su derecho a espacios públicos seguros y de calidad, los cuales vayan más allá de simples equipamientos o áreas verdes. En la literatura, encontramos un importante número de casos, los cuales proporcionan evidencia sobre la importancia que tienen los espacios públicos en el fomento de conductas sociales íntegras. Mulgan et. al. (2006) ponen de manifiesto la capacidad que tienen los espacios públicos en la generación de percepciones de seguridad, sentimiento de comunidad y confianza mutua. Chanes y Sanz (2012) afirman que los espacios públicos de calidad facilitan la adherencia a las normas, valores y códigos de conducta; los cuales permiten a los habitantes de un barrio o colonia una coexistencia armoniosa. Lo anterior queda en evidencia, cuando un conjunto de actividades matemáticas son realizadas en espacios públicos. Además de contribuir al desarrollo de las habilidades de pensamiento, este tipo de actividades contribuyen a la coexistencia vecinal armoniosa. En las investigaciones conducidas por Davies (2000), se demuestra que los espacios públicos correctamente diseñados, y que se encuentran bajo estrictos esquemas de mantenimiento, fomentan valores como la sensación de seguridad, la cohesión social y la satisfacción vecinal.

Por el contrario, autores como Mulgan et. al. (2006) y Lyndhurst (2004) demuestran que espacios públicos en condiciones de abandono se convierten en inductores para conductas violentas y antisociales. La teoría de las "*Ventanas Rotas*" propuesta por Wilson y Kelling (1982) sostiene que la apariencia física de los objetos y espacios públicos, influye en buena parte el comportamiento de los individuos. Por ejemplo, una casa con las "ventanas rotas" podría causar una percepción de abandono y por consiguiente, tener una mayor probabilidad de ser objeto de robos o saqueos. Autores como Chanes y Sanz (2012) y De Roda y Moreno (2001) han proporcionado evidencia empírica para afirmar que la permanencia residencial está directamente relacionada con el sentido de pertenencia vecinal. Es decir, la intención de residir en un determinado barrio está en función del grado de identificación del individuo con los valores y hábitos que ahí se preservan. Los espacios públicos representan los principales medios de contacto con otros vecinos, en los cuales se da un intercambio de información (e.g. consejos sobre uso del transporte público, recetas de cocina o escuelas de los hijos). De ahí la importancia que tienen los espacios públicos de calidad, para generar mayores niveles de confianza vecinal y fortalecer el tejido social.

Con base en lo anterior, es incorrecto afirmar que si dotamos a las áreas urbanas con espacios públicos amplios, limpios y funcionales; entonces los índices de violencia disminuirán automáticamente. Tampoco es cierto que los espacios públicos por sí solos, producirán un aumento de la interacción social y la comunicación vecinal. Ferguson y Mindel (2007), Chanes y Sanz (2012) y Bentley (1985) señalan que actividades que estimulen la apropiación del espacio público, favorezcan el apoyo vecinal y aumenten el nivel de percepción de seguridad, entre otras, son necesarias. Sobre esta premisa, se pone en evidencia la principal limitante del "Programa de Rescate de Espacios Públicos[1]", inicialmente a cargo de la Secretaría de Desarrollo Social en el año 2010, y actualmente bajo la responsabilidad Secretaría de Desarrollo Agrario, Territorial y Urbano. En los últimos dos años, el programa asignó el 16.3% de su presupuesto total para la realización de actividades sociales, lo cual es insuficiente. Por ejemplo durante los años 2014 y 2016, este programa organizó 890 eventos deportivos, 4,000 funciones de cine y 1,362 partidos de futbol. Esto representa un total de 6,252 actividades realizadas en espacios públicos para un país de aproximadamente 115 millones de habitantes. Es decir, 5.4 actividades sociales por cada cien mil habitantes en un periodo de dos años.

El objetivo de este trabajo consiste en investigar el impacto que tienen un conjunto de actividades matemáticas realizadas de forma sistemática en espacios públicos de colonias con marginación alta, en el fortalecimiento del tejido social. Hemos tomado como caso de estudio, la plaza principal de la colonia Chamilpa en Cuernavaca, México. La ciudad seleccionada se caracteriza por tener una precaria infraestructura de recreación y ocio, adicionalmente, la colonia Chamilpa se encuentra entre las de mayor marginación de la ciudad. Este trabajo intenta hacer una contribución original, en el sentido de aportar elementos adicionales a la investigación de los espacios públicos. Hasta donde sabemos, no existen investigaciones dentro del contexto Latinoamericano, las cuales evalúen los impactos de actividades matemáticas, en el bienestar psicológico y social de los usuarios de espacios públicos. El primer objetivo específico consiste en realizar una descripción sociodemográfica de los asistentes al evento. La evaluación en relación a si un evento de matemáticas contribuye a mejorar la cohesión social y la percepción de seguridad, se dejan como segundo objetivo. Las siguientes preguntas de investigación abordan tales objetivos: ¿Cuál es el perfil sociodemográfico de los asistentes a eventos matemáticos realizados en colonias urbanas con alta marginación? ¿Qué tipo de impacto tiene, tanto en la cohesión social como en la percepción de seguridad, la realización de actividades matemáticas en espacios públicos? La hipótesis en relación al segundo objetivo (pregunta) de investigación es: La utilización de los espacios públicos para realizar actividades matemáticas contribuye significativamente en el mejoramiento de la cohesión social, así como en la percepción de seguridad. Finalmente, sobre la base modelos estadísticos multivariados de reducción de dimensionalidad, evaluamos la relación existente entre la auto-percepción de riesgo o violencia, por una parte, y la organización-unión vecinal por la otra.

---

[1] Información detallada sobre este programa puede ser consultada en: http://www.gob.mx/sedatu/acciones-y-programas/programa-de-rescate-de-espacios-publicos

## 2. METODOLOGÍA

*2.1 Proyecto ARTEMAT: "Matemáticas para la Paz"*

Este es un proyecto innovador enfocado en fortalecer el tejido social y ampliar la visión que se tiene de las matemáticas. Esta iniciativa consiste en la instalación de estaciones de trabajo, en las cuales se desarrollan talleres, demostraciones, juegos y exposiciones. Todas las actividades presentadas están relacionadas con las matemáticas, dentro de un ambiente que estimula la creatividad y la convivencia. Uno de los objetivos específicos de ARTEMAT consiste en enseñar a los participantes, de forma lúdica e interactiva, el lado artístico de las matemáticas y a la vez, favorecer la interacción entre los vecinos de la colonia. El proyecto es financiado por el Programa Nacional de Prevención de la Violencia y el Delito (PRONAPRED, 2016), implementado por el Instituto de Matemáticas de la Universidad Nacional Autónoma de México y la organización civil Arte Sustentable (ver: http://artesustentable.org). En esencia, se buscaba mejorar la cohesión social y la visión que se tiene de las matemáticas.

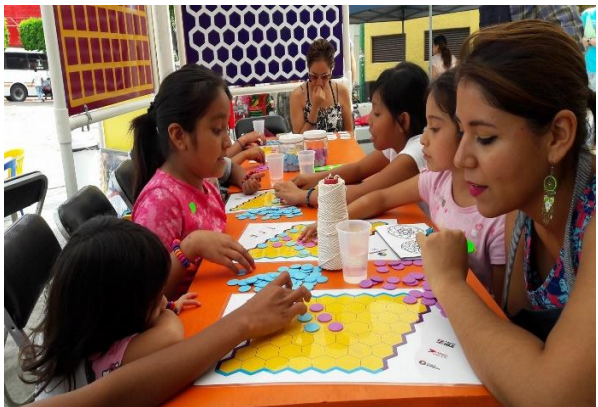 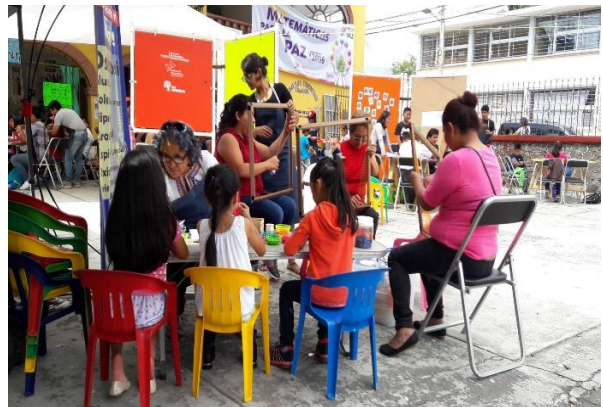

*Figura 1A*. Los asistentes.          *Figura 1B*. Apariencia del evento.

Desde la perspectiva de los participantes, ARTEMAT busca la creación de experiencias elocuentes, en donde cada asistente tenga un tiempo valioso, formativo y agradable, al mismo tiempo en que se estimula la creatividad e imaginación. El común denominador es que todas las actividades tienen un fundamento matemático, sobre el cual se busca mejorar el razonamiento lógico y de toma de decisiones. Entre las disciplinas matemáticas sobre las cuales se fundamenta ARTEMAT se encuentran la Geometría, Enigmas Matemáticos, Simetrías, Series, Lógica, Teselaciones, Teoría de gráficas, Topología, Fractales, Polígonos y Poliedros. En lo que respecta al diseño de las actividades, éstas se dirigen al público en general, con especial énfasis en niños y adolescentes. Las disciplinas matemáticas antes mencionadas, son enmarcadas sobre artes visuales y escénicas. Esta mezcla nos permite generar experiencias ricas, originales y atractivas entre los asistentes, las cuales a su vez tienen un efecto bidireccional. Por una parte, se busca un cambio de mentalidad hacia las matemáticas. Por la otra, un aumento en la percepción de seguridad entre los asistentes, así como un incremento en la cohesión social entre los vecinos (Ver Figuras 1A y 1B).

ARTEMAT tuvo lugar en la comunidad de Chamilpa, la cual se ubica al norte del municipio de Cuernavaca. Las razones para seleccionar esta comunidad (y no otras) obedecen a criterios de marginación y violencia. De acuerdo a datos del Diagnóstico Integral de Morelos, en el año 2016 esta comunidad se encontró en el polígono 1. Este último identificado como el más violento e inseguro, el cual reportó la tasa más alta violencia del municipio en 2015, igual a 8.4 delitos por cada mil habitantes (ver: http://morelosterritoriodepaz.org.mx). Por otra parte, los talleres se llevaron a cabo en seis estaciones de trabajo. La coordinación de cada estación estaba a cargo de dos personas con diferentes perfiles. Mientras un experto en matemáticas facilitaba el desarrollo de las habilidades del pensamiento analítico y numérico, un especialista en talleres comunitarios aseguraba que la actividad tuviera un alto contenido de creatividad e innovación. Previo al evento, los coordinadores recibieron una capacitación. Esto para asegurar que todos ellos estuvieran bien familiarizados con la misión y visión del proyecto, así como con el contenido de las actividades. Las actividades ofrecidas se presentan en la tabla I.

TABLA 1

Descripción de las actividades proporcionadas.

| Estación | Descripción |
| --- | --- |
| Pequeñitos | Actividades dirigidas a niños de 2 a 6 años. Se presentan temas como clasificación de objetos, pulseras y collares para tratar series, simetrías, figuras geométricas, tangrams[2]. |
| Hiloramas | En esta estación se dirige tanto a niños mayores de 8 años, así como a los adultos de cualquier edad. Aquí se presentan diferentes juegos con hilos para ilustrar conceptos matemáticos como tangentes, curvas, series, números primos y primos relativos. |
| Simetrías | Actividad enfocada a explicar el concepto de simetría (y sus implicaciones matemáticas). Para lo anterior, utilizamos caleidoscopios, mosaicos y espejos (Cunningham, et. al, 2015). |
| Enigmas, Magia y Nudos | Dinámicas grupales donde el razonamiento matemático es la clave para encontrar las soluciones. Aquí trabajamos principalmente con nudos (Eudave-Muñoz, 2006). |
| Poliedros y Polígonos | Las figuras geométricas regulares y sus propiedades matemáticas son aquí el tema central. Se realizan actividades colectivas, entre las cuales destacan pirámide de Sierpinski, papalote geométrico y las lámparas simétricas (Galarza y Seade, (2007). |
| Juegos | Esta estación se enfoca a realizar juegos de pareja o colectivo que involucren razonamiento matemático. Entre los juegos que aquí se ofrecen, se encuentran Nim, Hex y Sudokus (Galarza y Seade, (2007). |

---

[2] Consiste en un cuadrado compuesto de siete piezas, las cuales hay que ordenar para obtener figuras específicas (ver: http://www.disfrutalasmatematicas.com/definiciones/tangram.html).

El proyecto se llevó a cabo a lo largo de doce sesiones, comenzando el 21 de octubre y finalizando el 26 de noviembre de 2017. Cada sesión tuvo una duración de entre 6 y 8 horas. Con la finalidad de captar la mayor cantidad de asistentes, realizamos actividades de difusión previas. Se visitaron 16 escuelas primarias, 7 de ellas de turno matutino, 6 de turno vespertino y 3 más de tiempo completo, 2 escuelas secundarias de ambos turnos y 1 preparatoria de turno matutino. Más información sobre el proyecto ARTEMAT se encuentra en el Apéndice A.

*2.2 Instrumentos y Mediciones*

La Escala de Percepción del Control de La Violencia, desarrollada por Ramos (1990) mide las creencias que tiene un individuo, con respecto a quien tiene el control en hechos y/o sucesos delictivos. Los 21 items que componen la escala están dados en un intervalo tipo Likert con cinco valores posibles. Mientras un valor igual a uno hace referencia a "*En total desacuerdo*" el valor cinco es "*En total acuerdo*". La escala no tiene preguntas invertidas. El Proceso Analítico Jerárquico (PAJ), el cual consiste en realizar comparaciones por pares, aplicando los principios del álgebra lineal (Satty, 1980; Gómez y Cabrera, 2008), fue utilizado para priorizar los 21 items en base a su grado de contribución a nuestros objetivos de investigación. De acuerdo al algoritmo PJA, 8 items tuvieron una contribución alta, y por lo tanto, fueron incorporados a nuestro cuestionario. Adicionalmente, se incluyeron 2 items enfocados en medir el grado en que los asistentes perciben que las actividades matemáticas contribuyen, tanto al mejoramiento de la percepción de seguridad, como a la integración vecinal. Con lo anterior, nuestro cuestionario se compone de 10 items tipo Likert, los cuales se complementaron con 5 preguntas categóricas y referentes al perfil sociodemográfico de los asistentes.

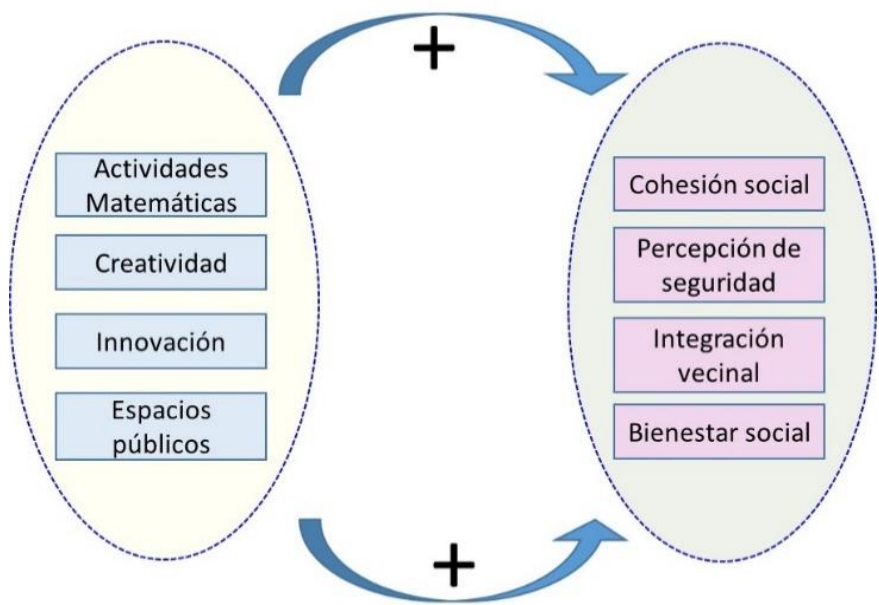

*Figura 2*. Mapa conceptual de la investigación.

Por otra parte, como resultado de nuestra revisión de literatura, identificamos formas en que actividades matemáticas, artísticas y culturales, realizadas en espacios públicos, impactan la integración vecinal y la cohesión social. Al respecto, otras investigaciones (Mallart y Deulofeu 2017; Arteaga, 2010; Mallart, 2008) ponen de manifiesto la importancia del pensamiento creativo y las artes, para lograr un alto impacto de las matemáticas durante festivales en espacios públicos. Sobre la base de este marco de referencia, construimos el mapa conceptual para este trabajo, el cual se encuentra alineado con la hipótesis de investigación (ver Figura 2).

*2.3 Métodos de análisis*

El análisis de los datos se realizó en dos fases. En la primera, realizamos un análisis descriptivo de los asistentes, obtuvimos las medidas de confiabilidad y validez del instrumento de medición, así como una exploración inicial de las relaciones entre las variables. Valores para los coeficientes de confiabilidad α de Cronbach y λ de Guttman (Cronbach, 1951; Tavakol y Dennick, 2011; Guttman, 1945), así como el Coeficiente de Correlación Inter-clase, se reportan como medidas de confiabilidad y validez del cuestionario. Un análisis de Componentes Principales (ACP) nos permite investigar las relaciones entre las variables estudiadas. De acuerdo con Smith (2002) el ACP es un método estadístico apropiado para investigar conjuntos de datos compuestos por un gran número de variables. En Barahona (2016) se menciona que una ventaja de ese método, sobre otras técnicas estadísticas multivariadas, es su capacidad para representar conjuntos de datos de gran tamaño, en visualizaciones intuitivas y relativamente fáciles de entender. A grandes rasgos, el ACP funciona de la siguiente forma. Sea $X$=[$X_1$, … ,$X_p$], una matriz con $p$ valores observables. A través de combinaciones lineales, calculamos aquellos componentes no-correlacionados, los cuales son capaces de capturar la mayor cantidad de variabilidad del conjunto inicial de datos *X*. De esta forma, los $p$ primeros componentes están dados por la siguiente formulación.

$$Y_1 = \mathbf{X}t_1, \quad Y_2 = \mathbf{X}t_2, \quad \ldots, Y_p = \mathbf{X}t_p \tag{1}$$

En (1), se asume que $var(Y_i)$ es valor máximo de varianza, dado $t_i' t_i = 1$, donde i=1,…,$p$. La covarianza entre cualquier par de componentes es $cov(Y_i, Y_j) = 0$, para i,j=1,…,$p$. De lo anterior, se dice que $Y_i$ y $Y_j$ son ortogonales. Si la variabilidad explicada en cada componente se expresa como $var(Y_i) = \lambda_i$, la variabilidad total del modelo es igual a $\mathbf{tr(S)} = \sum_i^p \lambda_i$, para i=1, 2, …, $p$. En relación al número total de componentes que deben ser incluidos, el criterio de Kaiser establece que aquel componente $Y_i$ con valor propio menor a 1.0 debería ser excluido (Kaiser, 1970). De esta forma, nuestro modelo se compone de los $m$ primeros componentes que satisfacen $\lambda_m > 1$; donde $\lambda_1 > \cdots > \lambda_m$ son los valores propios de la matriz de correlaciones **R**.

$$P = 100 \frac{\lambda_1 + \cdots + \lambda_m}{\lambda_1 + \cdots + \lambda_p} \tag{2}$$

En la formulación (2), *P* representa el porcentaje de varianza explicada, por aquellos componentes que satisfacen el criterio de Kaiser. De acuerdo con Greenacre y Hastie (1987), se espera que los primeros dos componentes capturen el 50% o más de la varianza total. Posteriormente, el índice de Tucker-Lewis y el Error Cuadrado Medio de los Residuales son proporcionados como medidas de confiablidad del ACP (McDonald y Marsh, 1990). Durante esta primera fase, los cálculos fueron realizados con el paquete estadístico R, versión 3.3.3 (ver: https://cran.r-project.org/).

Para la segunda fase del análisis se construyó un modelo de ecuaciones estructurales sobre la base del algoritmo de los Cuadrados Mínimos Parciales (CMP). Esta técnica estadística, propuesta por Wold (1985) y por Lohmöller (1989) se enfoca en maximizar la varianza explicada de las variables dependientes del modelo, y estimar los parámetros del modelo, siguiendo el algoritmo de los CMP. Este al ser un método no paramétrico, no requiere que los datos tengan ajuste de normalidad o bien, de alguna otra distribución. De acuerdo con Hair, Ringle y Sarstedt (2013), estos modelos se dividen en dos: el interno y el externo. Mientras el interno cuantifica las relaciones entre los constructos, el externo indica las relaciones entre las variables observadas y los constructos. En su forma matricial, el modelo interno puede expresarse de la siguiente forma:

$$\rho = \mathbf{X}\rho + \varepsilon \quad (3)$$

En la formulación (3), $\rho$ representa el vector de variables latentes, $\mathbf{X}$ es la matriz de coeficientes y $\varepsilon$ se asocia al error en el modelo. Por otra parte, a continuación aparece la expresión para el modelo externo.

$$y = \mathbf{W}y + \varepsilon \quad (4)$$

En (4), *y* denota el vector de cargas de las variables observables asociadas al constructo. $\mathbf{W}$ y $\varepsilon$ son la matriz de coeficientes y el error aleatorio, respectivamente. Nótese que para las ecuaciones (3) y (4) suponemos que los datos están estandarizados y por lo tanto se omiten los parámetros de intersección.

Los datos obtenidos en nuestra encuesta se utilizaron para obtener una solución numérica a las ecuaciones (3) y (4). Adicionalmente, el índice Kaiser-Meyer-Olkin (KMO), el Índice Normalizado de Ajuste (Bentler y Bonett, 1980) y el Error Cuadrático Medio Estandarizado, se proporcionan como medidas de ajuste y calidad del modelo. En la segunda fase se utilizó el paquete SmartPLS para realizar los cálculos (ver: https://www.smartpls.com/).

## 3. RESULTADOS

### *3.1 Perfil sociodemográfico de los asistentes*

En cuanto a las características sociodemográficas, podemos mencionar que el rango de edad de los asistentes estuvo entre los 3 y 65 años. En total recibimos a 1,302 personas, de los cuales el 29% tenía entre 3 y 11 años de edad. Un 54% entre un rango de 12 a 17 años. El 10% de los asistentes estaba entre 18 y 29 años. Finalmente, el 7% de ellos contaba con 30 años o más. En cuanto al género, 54% fueron mujeres y un 46% hombres. Adicionalmente, cada participante fue consultado si estaría interesado en responder nuestra encuesta. Durante las doce jornadas que llevamos a cabo, un total de 583 personas aceptaron a responder nuestra encuesta, de los cuales 337 fueron mujeres (58%), y 246 hombres (42%). De esta forma, recibimos en promedio 46 encuestas por día, con desviación estándar igual a ±8.5 (ver Figura 3).

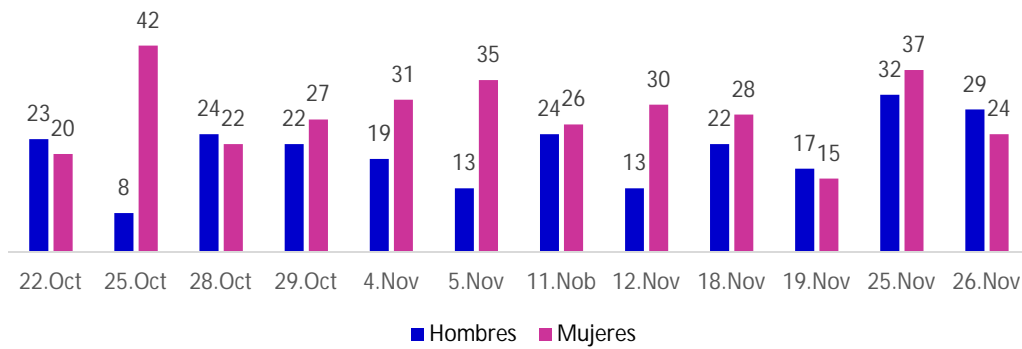

*Figura 3*. Encuestados por género.

Uno de los principales objetivos de nuestra encuesta consiste en identificar cambios en la percepción de los asistentes hacia las matemáticas, seguridad individual e integración vecinal. Por lo tanto, consultamos a todos los encuestados si estarían interesados en responder nuestro instrumento en más de una ocasión. Es decir, aplicaciones repetidas de un mismo cuestionario, pero en diferentes contextos y tiempos, nos permiten identificar el grado en que las actividades de ARTEMAT afectaban la percepción de los asistentes. Esto en la medida en que las semanas transcurrían. Aunque los cuestionarios se administraron anónimamente, evitando recopilar los nombres y domicilios de los asistentes, se les incluyó un ítem dicotómico, en la cual preguntamos si estaban contestando por primera vez la encuesta (SI o NO). Mientras que en la primera sesión de ARTEMAT el 100% de los encuestados respondieron que SI en esta pregunta, para la octava, este número disminuyó a 50%. De esta manera, en promedio el 30% respondió la encuesta en más de una ocasión, en contraste con el 70% que lo hizo solamente una vez (ver Tabla II).

TABLA II

Porcentaje de participantes que respondieron una sola vez (SI) el cuestionario.

|    | 22.Oct | 25.Oct | 28.Oct | 29.Oct | 4.Nov | 5.Nov | 11.Nov | 12.Nov | 18.Nov | 19.Nov | 25.Nov | 26.Nov |
|----|--------|--------|--------|--------|-------|-------|--------|--------|--------|--------|--------|--------|
| **SI** | 100%   | 96%    | 78%    | 84%    | 80%   | 65%   | 66%    | 49%    | 50%    | 66%    | 68%    | 55%    |
| **NO** |        | 4%     | 22%    | 16%    | 20%   | 35%   | 34%    | 51%    | 50%    | 34%    | 32%    | 45%    |

*3.2 Las matemáticas y el fortalecimiento de la cohesión social*

Previo al análisis multidimensional, revisamos las medidas de confiabilidad y estabilidad del instrumento. Los valores para los coeficientes α de Cronbach y λ de Guttman fueron iguales a 0.73 y 0.75, respectivamente. Así también, obtuvimos un Coeficiente de Correlación Interclase (CCI) igual a 0.72. Sobre la base de medidas de confiabilidad y estabilidad aceptables para el instrumento, procedimos a realizar el análisis factorial. El valor para el índice Kaiser-Meye-Olkin (KMO), el cual es un indicador del grado de adecuación de los datos para el análisis factorial, fue igual a 0.77. De acuerdo con Kaiser (1970) y Kaiser, (1974), valores superiores a 0.70 en el KMO indican buena adecuación. El método de extracción aplicado fue el de Análisis de Componentes Principales (ACP) con Varimax. Esta última es un tipo de rotación ortogonal, la cual maximiza la suma de las varianzas de la matriz cuadrada de cargas. El valor propio del primer componente fue de $\lambda_1=2.07$, mientras que el segundo tuvo un valor propio igual a $\lambda_2=1.46$. El índice de Tucker Lewis, el cual es una medida de bondad de ajuste para el ACP fue igual a 0.89, así también, el Error Cuadrático Medio resultó igual a 0.06. Los diez ítems incluidos en el estudio mostraron cargas de 0.41 y superiores en los dos componentes retenidos. Así también, la matriz de carga muestra las propiedades convergencia y discriminación (ver Tabla III).

TABLA III

Análisis de Componentes Principales (matriz de cargas).

| item | PC1 | PC2 | Comunalidades |
|------|-----|-----|---------------|
| X8_ActiHoyEviDelitos | 0.676 |       | 0.48 |
| X7_ContribResolSeg   | 0.630 |       | 0.40 |
| X6_OrgaVecSeg        | 0.614 |       | 0.38 |
| X10_ActPubcColSegur  | 0.531 |       | 0.31 |
| X2_UniVecViole       | 0.493 |       | 0.25 |
| X1_LlevBienMolest    | 0.431 |       | 0.20 |
| X9_RelVecDelito      |       | 0.566 | 0.37 |
| X4_MalaRelVec        |       | 0.621 | 0.40 |
| X5_SimpSegurPer      |       | 0.727 | 0.53 |
| X3_PoliciaSeg        |       | 0.405 | 0.17 |
| Valores Propios      | 2.07  | 1.46  |      |

Posteriormente se calcularon las correlaciones de los ítems con los dos primeros componentes, los cuales se representaron en un círculo de correlaciones (ver Figura 4). Al ser esta una representación vectorial, existe una interpretación sobre la magnitud, dirección y sentido de los ítems. En el primer cuadrante encontramos ítems que tienen en común los temas de la seguridad y la violencia. Por lo tanto, los hemos agrupado en una variable latente denominada "*Auto-percepción de seguridad*", la cual responde a la siguiente la pregunta *¿Qué tan seguro me siento?* Por otra parte, en el cuarto cuadrante encontramos ítems relacionados con la cohesión e integración vecinal. Esta segunda variable latente la denominamos "*Integración vecinal*" y responde a la pregunta: *¿Qué tan integrado estoy con mis vecinos*? La dirección de los vectores puede ser interpretado así: entre más a la derecha encontremos posicionado el vector, indica mayor intensidad de tal variable. Por ejemplo, el ítem "*X10. Las actividades matemáticas en lugares públicos pueden mejorar la seguridad en mi colonia*" es la más cercana a la parte derecha del primer componente. Por el contrario, la pregunta "*X3 La policía es la única que puede darme seguridad*" resulta la más alejada al primer componente.

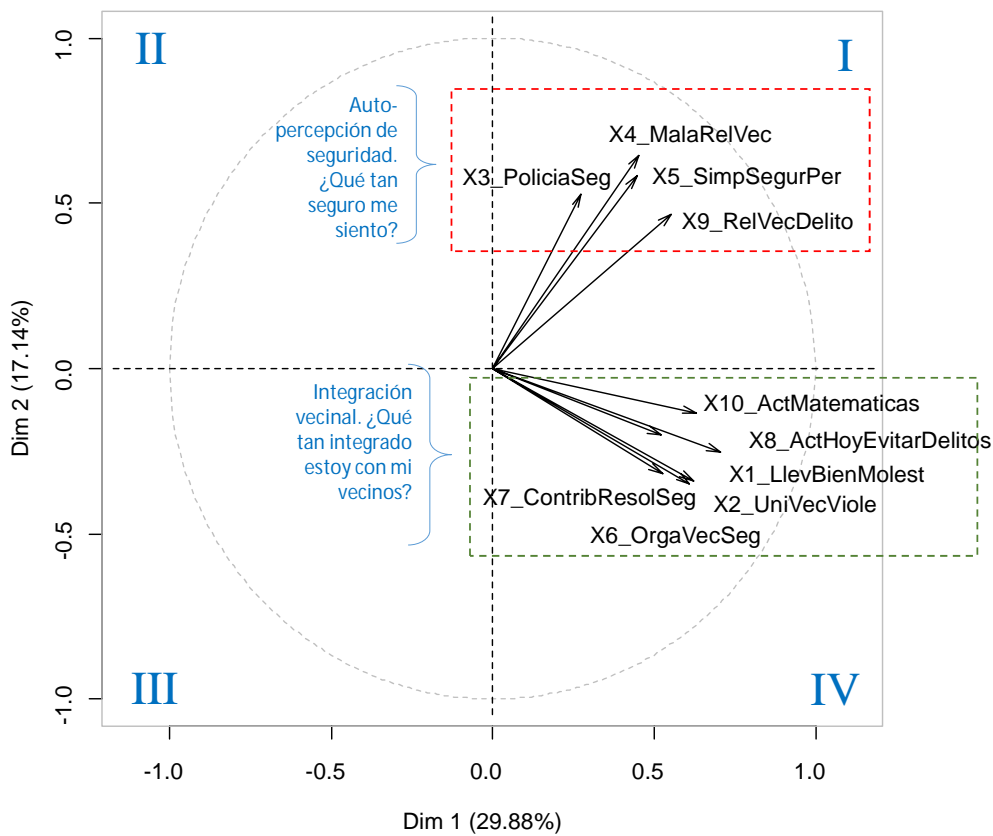

*Figura 4.* Círculo de correlaciones en su forma vectorial.

Posteriormente, obtuvimos una nube de puntos sobre el plano factorial, la cual está dada por la contribución que tiene cada uno de los individuos (n=583) a los componentes extraídos. La variable fecha se incorpora al análisis como variable suplementaria. Esta, aunque tiene una coordenada del tipo c($x,y$) en el plano, al ser suplementaria no contribuye a la generación de los valores propios de los ejes (componentes). Lo anterior nos permite crear visualizaciones de las fechas con respecto a los individuos. Por ejemplo, dos fechas aparecerán muy cercanas en el plano, si obtuvieron resultados semejantes en la encuesta. Por el contrario, dos fechas estarán lejanas, cuando sus resultados en la encuesta sean diferentes. En la Figura 5, aparecen las doce fechas, en las cuales se llevó a cabo el evento ARTEMAT. Por una parte, en el tercer cuadrante (III) encontramos las fechas iniciales del evento (22, 25 y 28 de octubre). Por la otra, las fechas para las últimas jornadas aparecen en los cuadrantes primero (I) y cuarto (IV). Estas son 12, 19 y 25 de noviembre.

Una interpretación conjunta de los gráficos mostrados en las Figuras 4 y 5 se realiza de la siguiente forma: mientras en el primer cuadrante (I) encontramos los ítems relacionados con la autopercepción de seguridad, las fechas que aparecen en este mismo cuadrante son el 29 oct, 25-nov y 26-nov. Por otra parte, los ítems concernientes a la integración y cohesión vecinal se encuentran en el cuarto cuadrante, y los cuales tienen de contraparte las jornadas 12 y 19 de noviembre. De esta forma, se muestra que las fechas se han movido de izquierda a derecha sobre el eje horizontal. Con lo anterior, se proporciona evidencia estadística parcial, la cual soporta que la realización de actividades matemáticas en espacios públicos contribuye al mejoramiento de la cohesión social y la integración vecinal.

Así también, en la Figura 5 aparece una trayectoria, la cual es representada por flechas rojas que se encuentran entre los días. Las puntas de las flechas nos indican la dirección de la trayectoria. Mientras al inicio encontramos las primeras jornadas para ARTEMAT (flechas "a" y "b"), en las últimas aparecen las fechas finales ("e" y "f"). Nótese que el plano horizontal está caracterizado principalmente por las preguntas "*X8. Actividades como la de hoy evitan que sucedan actos delictivos*" y "*X10. Las actividades matemáticas en lugares públicos pueden mejorar la seguridad en mi colonia*"; además, la trayectoria de las fechas rojas es de izquierda a derecha sobre el mismo plano horizontal. Es decir, las fechas iniciales (22 y 25 Oct) se encuentran más alejadas de los vectores para las preguntas *X8* y *X10*. En contraste, las fechas finales del evento (25 y 26 de Nov) se encuentran más cercanas a tales vectores (*X8* y *X10*). Lo anterior nos proporciona evidencia adicional para afirmar que ARTEMAT tuvo un impacto en la cohesión social y la integración vecinal. En la medida en que fueron transcurriendo las jornadas del evento, la percepción de los individuos, sobre si la realización de actividades matemáticas en lugares públicos contribuye a la cohesión social, fue cambiando.

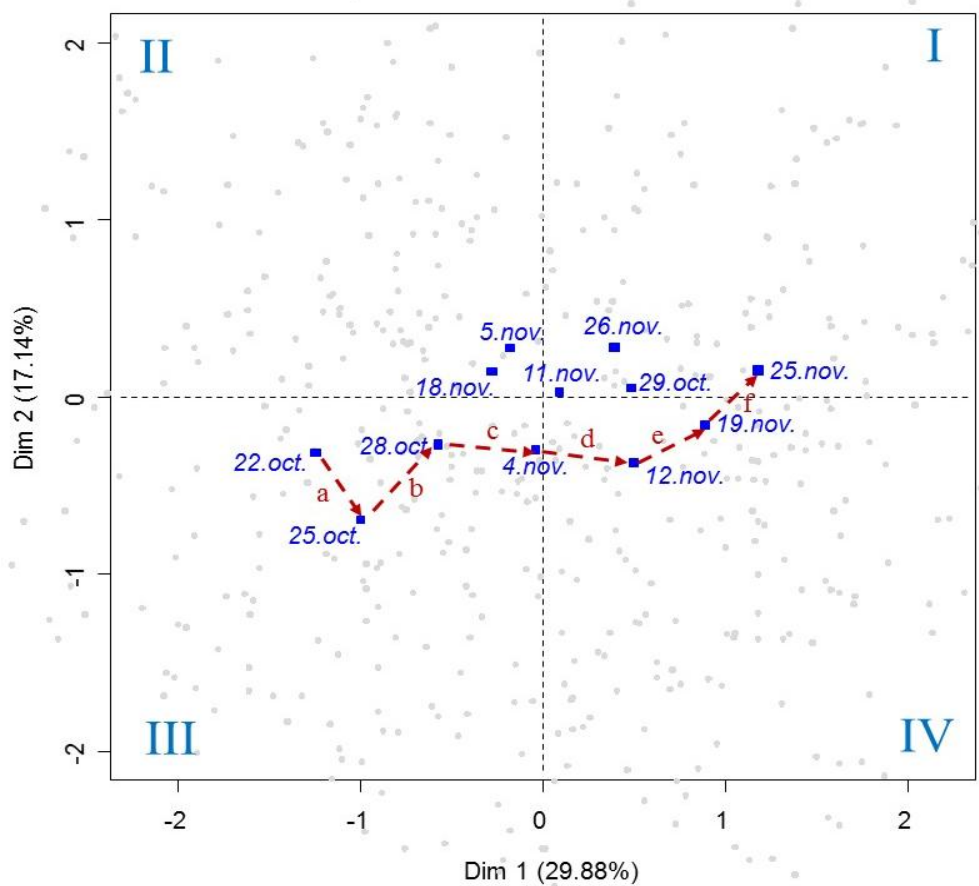

*Figura 5.* Nube de puntos y fechas sobre el plano factorial.

La última parte de esta sección se enfoca en los resultados del modelo de ecuaciones estructurales. Para tales efectos, agrupamos a los ítems X8 y X10 en una variable latente denominada "*Actividades matemáticas en espacios públicos*". Por otra parte, tres ítems "*X2 La violencia puede disminuir si existe unión entre los vecinos de mi colonia*", "*X6 Si los vecinos de mi colonia nos organizamos, podríamos evitar que sucedieran actos delictivos*" y "*X7 Si me lo propongo, puedo contribuir en algo para resolver el problema de la violencia*" conformaron la variable latente denominada "*seguridad y unión vecinal*". Las razones para incluir estos ítems y no otros en nuestro modelo estructural, se basan en la definición operativa de variables realizada para tales efectos. Desde una perspectiva conceptual, nótese que los cinco ítems incluidos en el modelo hacen referencia a constructos como actividades matemáticas, unión vecinal, violencia, seguridad y actos delictivos. Sobre este marco, el objetivo es arrogar evidencia cuantitativa que soporte o refute que las actividades matemáticas mejoran la percepción de seguridad, así como, la unión vecinal.

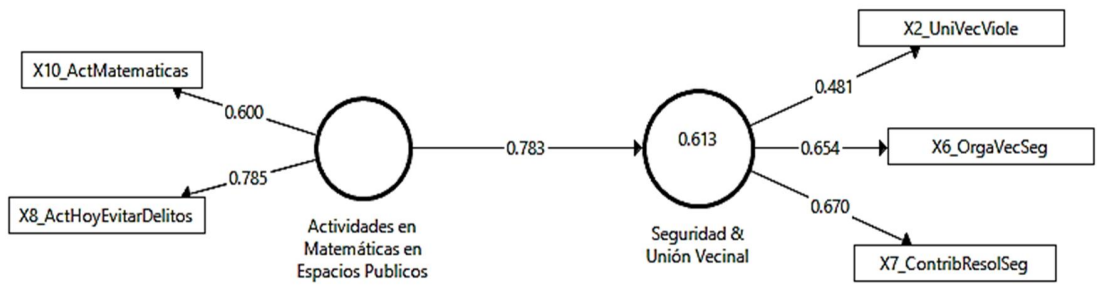

*Figura 6.* Modelo de ecuaciones estructurales.

En la figura 6 aparece un diagrama para el modelo de ecuaciones estructurales. En este caso nuestra variable independiente son "*las actividades matemáticas en espacios públicos*". Por otra parte, la variable de respuesta es la "*seguridad y unión vecinal*". Como se puede apreciar, tales constructos se encuentran unidos por el coeficiente de trayectoria (path coefficient). Este a su vez indica la dirección de la relación (en este caso de izquierda a derecha). Valores cercanos a 1.0 indican una relación fuerte. Nuestro coeficiente es igual a 0.783. Lo anterior aporta elementos a favor de nuestra hipótesis inicial: las actividades matemáticas en espacios públicos mejoran la percepción de seguridad y unión vecinal.

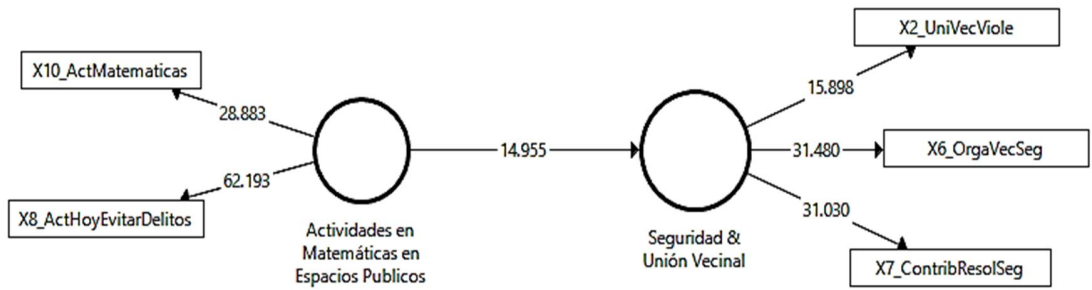

*Figura 7.* Pruebas de hipótesis para del modelo.

Adicionalmente, calculamos el estadístico t-student para todos los coeficientes que integran el modelo. Lo anterior para evaluar la hipótesis nula $H_o$: los coeficientes son significativamente diferentes de cero. Con una significancia de α=0.05, el umbral de rechazo es menor o igual a 1.96. Es decir, vamos a rechazar la hipótesis de significancia para aquellos coeficientes en donde el estadístico t-student sea menor a 1.96. En nuestro caso, fallamos en rechazar $H_o$ en todos los casos, y por tanto, los coeficientes en el modelo son significativos (ver figura 7). Por último, como medidas bondad de ajuste, se reportan el Índice Normalizado de Ajuste (INA), igual a 0.998 y el Error Cuadrático Medio Estandarizado con valor de 0.008.

## 4. DISCUSIÓN

En relación al perfil sociodemográfico de los asistentes, mediante nuestro análisis descriptivo identificamos que se encuentran en un rango de edad de 3 y 65 años. Además, 29% de los participantes estaba entre los 3 y 11 años de edad. La mayoría (54%) estaba entre un rango de 12 y 17 años de edad. El 10% estuvo entre los 18 y 29 años. De esta forma, identificamos que la población que asiste a este tipo de actividades matemáticas en parques y espacios públicos, son principalmente jóvenes menores de 30 años. En este caso, tal grupo representa el 93% del total de los visitantes. Lo anterior es congruente a lo observado durante otros eventos similares y relacionados con la difusión de las ciencias. Por ejemplo, durante la Semana Nacional de Ciencia y Tecnología, la cual es el evento más grande sobre la difusión de la ciencia que se realiza en México, se recibieron durante su edición del año 2016 alrededor de 120,000 visitantes, de los cuales 80% estaban en un rango de edad de entre los 3 y 30 años (CONACYT, 2016). En lo que respecta al género de los asistentes, encontramos una mayor participación de las mujeres (54%) con respecto a los hombres (46%). El hecho de que las mujeres muestren mayor entusiasmo a participar en eventos científicos y culturales es una tendencia internacional. Por ejemplo, en la Encuesta de Hábitos y Prácticas Culturales en España 2014-2015, realizada por el Ministerio de Educación, Cultura y Deporte, las mujeres aparecen por encima de los hombres en la asistencia a eventos culturales en espacios públicos (MECD, 2016). De igual forma, la Encuesta Nacional de Consumo Cultural de México, llevada a cabo en el año 2012 reafirma esta tendencia (INEGI, 2012). Del total de asistentes a eventos culturales en espacios públicos para ese año, un 54% fueron mujeres, mientras que el 46% resultaron hombres.

Otro hallazgo relevante es relación a la percepción de inseguridad existente entre los hombres y las mujeres asistentes al evento. Los resultados indican que los hombres se auto-perciben más inseguros que las mujeres. En la pregunta cuatro del cuestionario que dice: "*X4. Las cosas malas que me suceden dependen de cómo me llevo con los vecinos*" los hombres obtuvieron una media igual a $\mu_h$=2.83 ±1.4. En contraste, las mujeres que lograron una media de $\mu_m$=2.41 ± 1.38. Así también, para la pregunta cinco "*X5. Mi seguridad personal depende de qué tan simpático y agradable soy con los vecinos*", las medias fueron $\mu_h$=2.78 ±1.41 y $\mu_m$= 2.44 ± 1.32 para los hombres y mujeres respectivamente. Lo anterior es consistente con las cifras sobre homicidios reportadas por el Instituto Nacional de Estadística y Geografía en el 2015. Durante ese año se registraron 20,762 homicidios, de los cuales 88% correspondieron a hombres y el 12% a mujeres (INEGI, 2015). De esta forma, se pone de manifiesto que la violencia en México es principalmente un problema que afecta al género masculino. Ser hombre joven en México se asocia una mayor probabilidad de estar relacionado con algún tipo de suceso violento, ya sea como víctima o perpetrador.

Uno de los objetivos de investigación para este trabajo es identificar en qué medida las actividades matemáticas llevadas a cabo en espacios públicos contribuyen a mejorar la percepción de seguridad y la cohesión social. En este sentido, el análisis de componentes principales realizado a nuestros datos nos permitió definir dos constructos. El primero relacionado con la auto-percepción de seguridad que tienen los participantes en el evento. El segundo dado por la integración vecinal y las actividades matemáticas en espacios públicos. A través tal análisis, identificamos correlaciones positivas entre ambos grupos de variables. Es decir, a medida que aumentan la integración vecinal y las actividades matemáticas en los espacios públicos, se registra un incremento en los niveles de auto-percepción de seguridad. Aunque no identificamos otros trabajos que investiguen el impacto de los festivales matemáticos en la percepción de seguridad, nuestros resultados se pueden contrastar con otros similares. Al respecto, Albanese y Perales (2014) estudian la relación existente entre pensamiento matemático, artesanías y espacios públicos. Carro, Valera y Vidal (2010) proporcionan un marco conceptual, en el cual diferentes constructos afectan, ya sea de manera positiva o negativa, la percepción seguridad en los espacios públicos. Por ejemplo, la falta de iluminación de los espacios disminuye la percepción de seguridad. Por el contrario, la utilización de los mismos con fines lúdicos y culturales se asocia mayores niveles en la percepción de seguridad. Adicionalmente, Valera y Guàrdia (2014) validaron el marco teórico antes mencionado con datos obtenidos mediante una encuesta, y los cuales se analizaron con un modelo de ecuaciones estructurales. En tal estudio, se proporciona evidencia, la cual soporta que la percepción de inseguridad en espacios públicos es explicada en términos estructurales por tres variables latentes: las características físicas del propio espacio, el nivel de socialización los vecinos y las competencias propias del individuo.

Adicionalmente, en el estudio realizado por Ferraro (1995) se analiza el impacto que tiene el hecho de compartir información sobre un espacio público en la percepción de seguridad. Los resultados evidenciaron que la percepción de seguridad disminuye en la medida en que los individuos comparten información sobre ciertos hechos delictivos. Así también, en el trabajo realizado por Tejera (2012) se llevaron a cabo observaciones sistemáticas durante cuatro meses en seis diferentes espacios públicos. Esto con la finalidad de identificar diferentes niveles de percepción de seguridad. En sus resultados, el autor encontró diferencias significativas entre los parques estudiados, así como evitación del espacio público por parte de mujeres y personas mayores. Por otra parte, en nuestro análisis de componentes principales presentado en la Figura 5, igualmente se analizan cambios en el nivel de percepción de seguridad con respecto al tiempo. En este caso aparecen en el primer plano factorial, tanto la nube de individuos como todas las fechas en las que realizamos el evento matemático. Tal como se ilustra con las flechas identificadas con las primeras seis letras del alfabeto, los niveles de percepción en la seguridad así como la cohesión social aumentaron, conforme transcurría el evento. En la medida en que los asistentes socializaban su experiencia lúdica con las matemáticas, los niveles de percepción en la seguridad e integración vecinal aumentaban. Se aprecia que las jornadas iniciales del evento (22 y 25 de octubre) aparecen en el tercer cuadrante y más alejadas a las variables de auto-percepción de seguridad e integración vecinal. Por el contrario, las jornadas finales del evento (25 y 26 de noviembre) aparecen más cercanas a nuestras variables de interés.

Partiendo de la premisa que "*correlación no implica causalidad*"[3] el análisis de componentes principales discutido en párrafos anteriores, está limitado y proporciona evidencia parcial a favor de nuestra hipótesis de investigación. Lo anterior debido a que está basado en la matriz de correlaciones. Con la finalidad de sopesar esta limitante, se propone un modelo de ecuaciones estructurales. A diferencia de los análisis de correlaciones tradicionales, los modelos estructurales y sus coeficientes de camino nos permiten identificar causas y efectos. En la Figura 6 se presenta el diagrama para nuestro modelo estructural. Todos los ítems mostraron cargas iguales 0.481 y superiores en las viables latentes. El coeficiente de camino entre los constructos "*Actividades matemáticas en espacios públicos*" y "*Seguridad y unión vecinal*" fue igual a 0.783, con una dirección de derecha a izquierda. De acuerdo con Hair, Ringle y Sarstedt (2013) y con Hair et. al (2014), valores a 0.70 indican causalidad fuerte entre las variables estudiadas. Para hacer los resultados más robustos, realizamos una simulación son 500 submuestras a través del método "Bootstrapping"[4] (Preacher y Hayes, 2008; Efron, 1982). Con un nivel de significancia α=0.05, todos los coeficientes del modelo resultaron significativos. Estos resultados complementan los obtenidos en otras investigaciones. En Marzbali et. al (2012) se valida un modelo estructural, en el cual cuatro constructos tienen un impacto significativo en la seguridad de los espacios públicos. Ferraro (1995), Tejera (2012), y Valera y Guàrdia (2014) discuten la relación que existe entre los diferentes usos que tienen los espacios públicos y la cohesión social.

---

[3] "Con esto, por tanto a causa de esto" (En latín "*Cum hoc ergo propter hoc*") es una falacia que se comete al inferir mediante los coeficientes de correlación que dos o más eventos están conectados causalmente.

[4] Método de remuestreo empleado para aproximar un conjunto de datos a una distribución. También se usa para estimar el error y varianza de un modelo estadístico.

## 5. REFLEXIONES FINALES

Este trabajo investiga el impacto que tienen un conjunto de actividades matemáticas, en dos variables latentes. Mientras la primera se refiere a la integración vecinal, la segunda está relacionada con los niveles de autopercepción de seguridad. Los resultados aquí obtenidos hacen una contribución original, al aportar elementos adicionales a la discusión sobre el uso de los espacios públicos en Latinoamérica. Estos pueden ser de utilidad para investigaciones posteriores en el tema, las cuales se enfoquen en la formulación de marcos teóricos o metodológicos. O bien, para aquellos trabajos que busquen explicar sistemáticamente las asociaciones existentes entre matemáticas educativas y bienestar psicológico de usuarios de espacios públicos.

Por otra parte, los tomadores de decisiones y responsables de la formulación de políticas públicas encaminadas a mejorar la seguridad en espacios públicos urbanos, tienen en este trabajo elementos adicionales de análisis. El Programa de Rescate de Espacios Públicos representa actualmente la principal herramienta que tiene el gobierno de México para intervenir espacios públicos en el país (SEDATU, 2013). Al respecto, el informe elaborado por el Consejo Nacional de Evaluación de la Política de Desarrollo Social (CONEVAL) en relación al desempeño de este programa en el 2013, resalta la necesidad de involucrar a la propia comunidad a tener una participación más activa en las actividades realizadas en los espacios públicos (CONEVAL, 2013). Los resultados presentados en este trabajo aportan elementos de análisis, los cuales van en la misma dirección que las recomendaciones del CONEVAL. El nivel de impacto positivo de una actividad realizada en espacios públicos, estará en relación directa con el nivel de involucramiento de los asistentes. En el contexto latinoamericano los tomadores de decisiones intentan solucionar el problema de inseguridad en los espacios públicos bajo el paradigma tradicional (limitado) de enfocarse casi exclusivamente en mejorar la infraestructura. La mayoría de las políticas públicas parten de la premisa que dotar a los espacios públicos con mayor infraestructura física (cámaras, estaciones de vigilancia, puertas de acceso, etc), tendrá como resultado una mejora en la percepción de seguridad. Lo anterior es una visión fragmentada del problema. Al respecto, es necesario desarrollar más investigaciones sobre el tema, las cuales ayuden a mejorar la comprensión del impacto que tienen las actividades recreativas realizadas en espacios públicos, sobre la percepción de seguridad e integración vecinal. Un mejor entendimiento del impacto que tienen los diferentes usos de los espacios públicos, ayudará a tomar mejores decisiones sobre las políticas públicas adecuadas para su intervención y mejoramiento.

Por último, es importante mencionar que el alcance de nuestros resultados está limitado, tanto por el tamaño de la muestra, como por el alcance geográfico. Si bien los resultados aquí presentados únicamente caracterizan la situación actual de una colonia en particular en México, estos pueden servir de base para realizar comparaciones con otros estudios similares, o bien como punto de referencia para investigaciones posteriores.

# REFERENCIAS BIBLIOGRÁFICAS

# Apéndice A

ARTEMAT. Matemáticas para la Paz.

Memoria fotográfica.

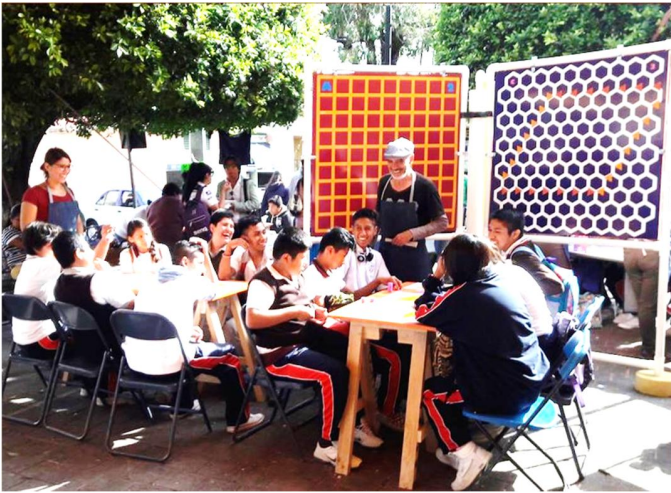

*Figura A1*. Apariencia de los asistentes en las mesas de trabajo

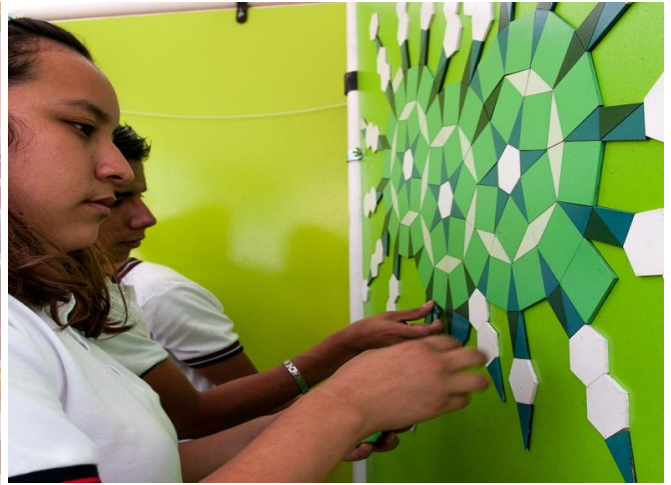

*Figura A2*. Apariencia de los asistentes desarrolando las actividades. Foto: Aubin Arroyo

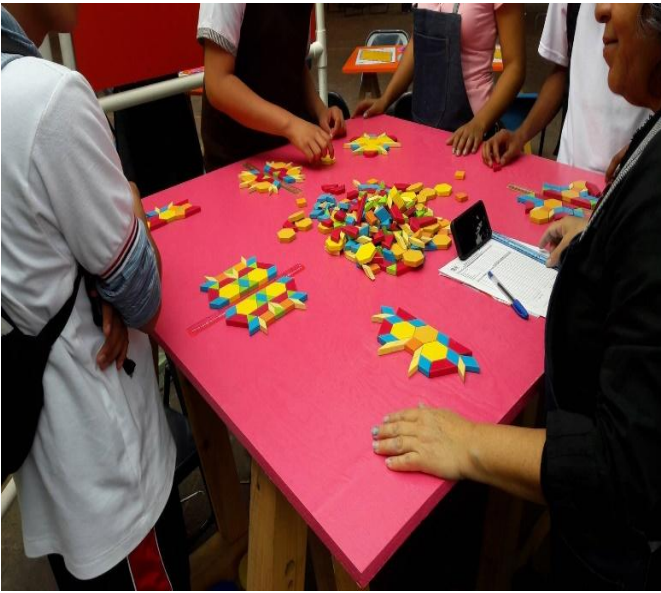

*Figura A3*. Apariencia de las mesas de trabajo

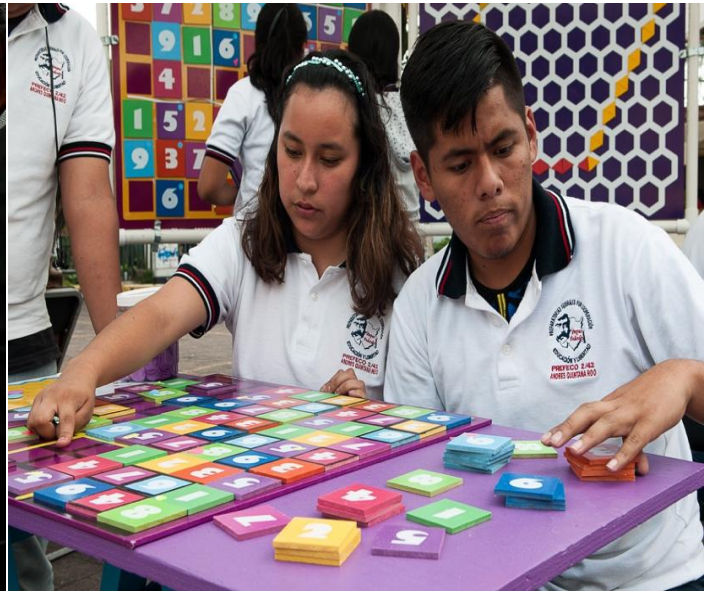

*Figura A4*. Apariencia de una de las actividades implementadas. Foto: Aubin Arroyo

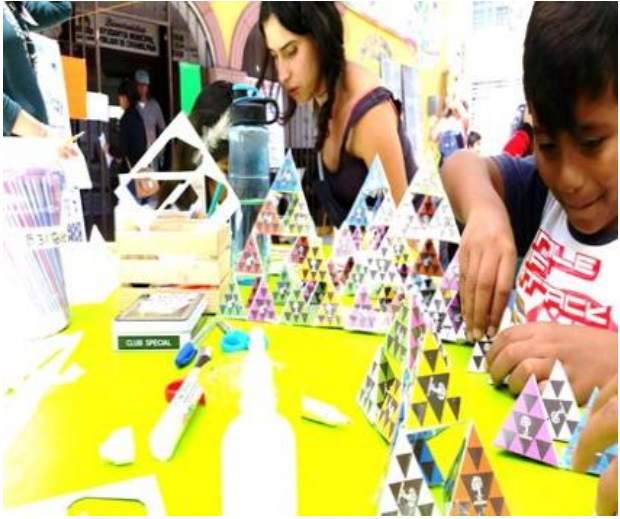
*Figura A5.* Apariencia de las actividades

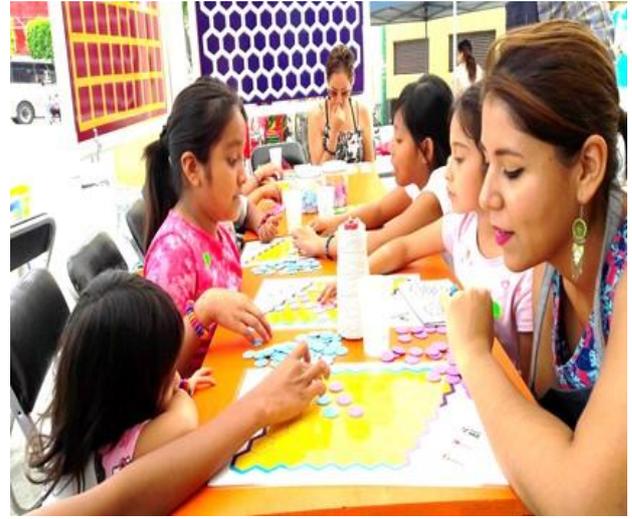
*Figura A6.* Apariencia de las actividades

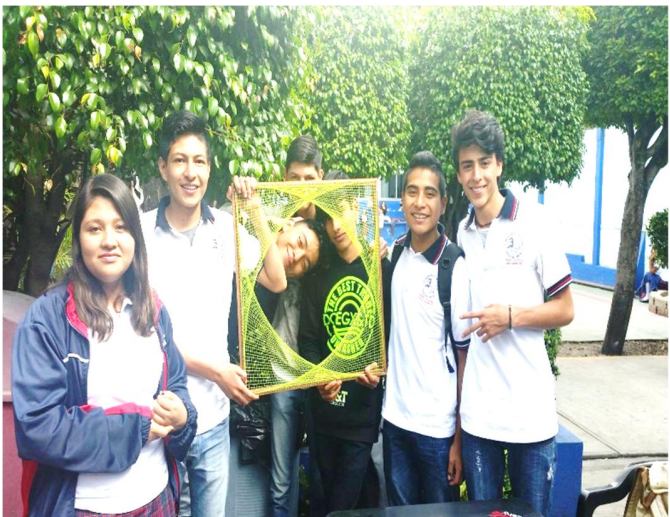
*Figura A7.* Los voluntarios de la preparatoria

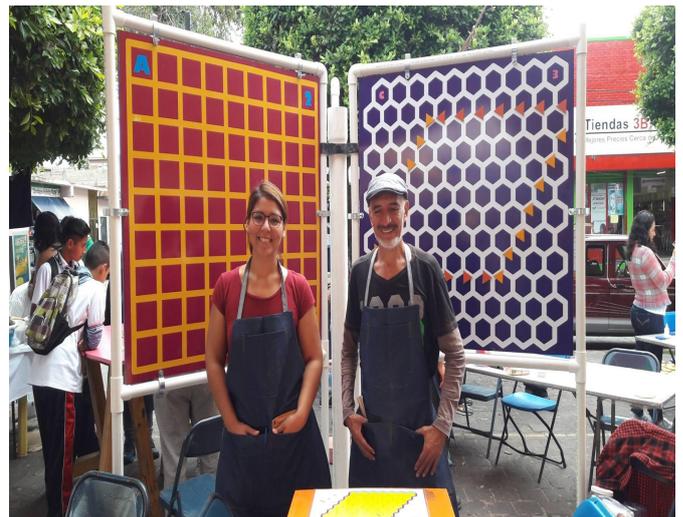
*Figura A8.* Los coordinadores de actividad (Talleristas)